\documentclass[reqno,12pt]{amsart}
\usepackage{mathrsfs, amssymb, latexsym, amsmath, graphicx, psfrag}
\newcounter{thm}

\newtheorem{theorem}[thm]{Theorem}
\newtheorem{corollary}{Corollary}
\newtheorem{lemma}[corollary]{Lemma}
\newtheorem{rmk}[corollary]{Remark}

\newcommand{\bk}{\mathrm{bk}}
\newcommand{\nzbk}{\mathrm{nzbk}}
\newcommand{\pa}[1]{\Pi^{A}_{#1}}
\newcommand{\ba}[1]{\Gamma^{A}_{#1}}
\newcommand{\pb}[1]{\Pi^{B}_{#1}}
\newcommand{\bb}[1]{\Gamma^{B}_{#1}}
\renewcommand{\a}{\alpha}
\renewcommand{\d}{\delta}
\newcommand{\bi}[1]{\varphi_{#1}}
\newcommand{\la}{\langle}
\newcommand{\ra}{\rangle}

\title[Temperley-Lieb algebra and chromatic joins]{The Gram matrix 
of a Temperley-Lieb algebra is similar to the
matrix of chromatic joins}
\author{Qi Chen\ \ and \ \ Jozef H. Przytycki}

\begin{document}
\maketitle

\centerline{Dedicated to the memory of Xiao-Song Lin (1957-2007)}

\vspace{10ex}

\centerline{{\em To appear in Communications in Contemporary Mathematics}}

\section{Introduction}\label{sec1}
\addtocounter{footnote}{1}

Rodica Simion noticed experimentally
that matrices of chromatic joins (introduced by
W. Tutte in \cite{Tu2}) and the Gram matrix of the Temperley-Lieb
algebra, have the same determinant, up to renormalization.
In the type $A$ case, she was able to prove this by comparing the known 
formulas: by Tutte and R. Dahab \cite{Tu2,Dah}, in the case of 
chromatic joins, and by P. Di Francesco, and
B. Westbury \cite{DiF,We} (based on the work by K. H. Ko and 
L. Smolinsky \cite{KS}) in the Temperley-Lieb case; see \cite{CSS}. 
She then asked for a direct 
proof of this fact \cite{CSS}, \cite{Sch}, Problem 7. 

The type $B$ analogue was an open problem central to the work of Simion
\cite{Sch}. She demonstrated strong evidence that the type $B$ 
Gram determinant of the Temperley-Lieb algebra is equal 
to the determinant 
of the matrix of type $B$ chromatic joins, after a substitution similar to 
that in  type $A$, cf. \cite{Sch}.

In this paper we show that the matrix $J_n$ of chromatic joins and 
the Gram matrix $G_n$ of the Temperley-Lieb algebra are similar 
(after rescaling), with the change of basis given 
by  diagonal matrices. More precisely we prove the following two results:

\begin{theorem}\label{A}
We have $J^A_n(\d^2)=PG^A_n(\d)P$, where $P=(p_{ij})$ is a diagonal matrix
with $p_{ii}(\d)= \d^{\bk(\pi_i)-n/2}$; here $\bk(\pi_i)$ denotes the 
number of blocks in the type $A$ non-crossing $n$-partition $\pi_i$; see 
2.1 of Section 2 for precise definitions. 
\end{theorem}

\begin{theorem}\label{B}
We have
$J^B_n(\d^2)=P^BG^B_n(1,\d)P^B$ where $P^B=(p^B_{ij})$
is a diagonal matrix with $p^B_{ii}(\d)= \d^{\nzbk(\pi_i)-n/2}$;
here $\nzbk(\pi_i)$ denotes half of the number of non-zero blocks 
in the type $B$ non-crossing $n$-partition $\pi_i$; see 
2.2 of Section 2 for precise definitions.
\end{theorem}

\section{Definitions and notation}
In the description of the matrix of chromatic joins (of type $A$ and $B$) 
we follow V. Reiner, R. Simion and F. Schmidt \cite{Rei,Sim,Sch}.

\subsection{The type $A$ case}
\begin{figure}[ht]
        \centering
                \includegraphics[height=2.5cm]{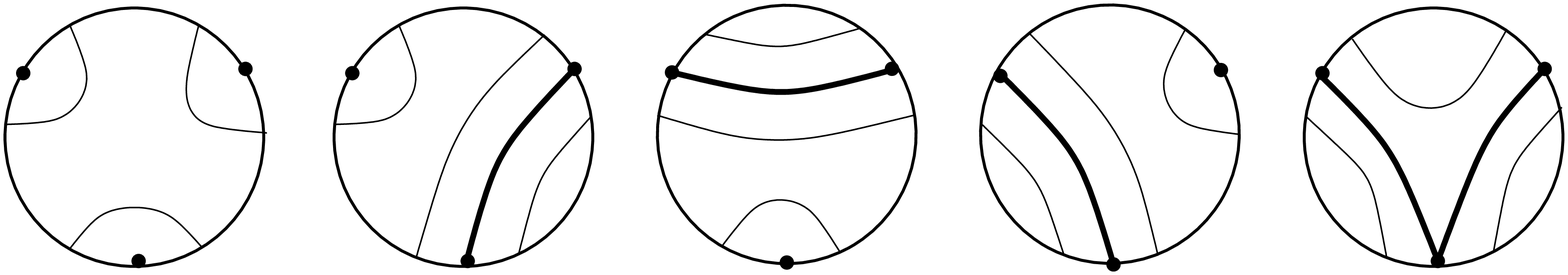}
\caption{The bijection $\bi A: \pa 3\to\ba 3$.}\label{BijectionA3}
\end{figure}

An {\em $n$-partition of type $A$} is a partition $\pi$ of the $n$ element 
set $\{1,2,...,n\}$ into blocks. The number of blocks is denoted by 
$\bk(\pi)$.
To represent $\pi$ pictorially, we place
the numbers $1,2,...,n$ anti-clockwise around the boundary circle 
of the unit disk and draw a chord, called a {\em connection chord}, 
in the disk between two numbers $i<j$ if they are
in the same block of $\pi$ and there is no $k$ in
the same block with $i<k<j$. We say that $\pi$ is
{\em non-crossing} if all connection chords can be drawn 
without crossing each other.  Notice that each block is represented by a tree. 
Denote the set of all non-crossing $n$-partitions of type $A$ 
by $\pa n$. On the other hand if $n=2m$ is even, we have  
{\it bipartitions} of $2m$ points of type $A$, those $2m$-partitions 
of type $A$ with every block 
containing exactly 2 numbers. 
Denote the set of all non-crossing $2m$-bipartitions of type $A$ by $\ba m$.
We have a bijection $\bi A : \pa n \to\ba n$ realized by considering 
the boundary arcs
of a regular neighborhood of the connection chords (see Fig. \ref{BijectionA3} 
and compare Fig. \ref{noncrB2}, Fig. \ref{figBbijection}).

\subsection{The type $B$ case}
\begin{figure}[ht]
        \centering
                \includegraphics[height=2cm]{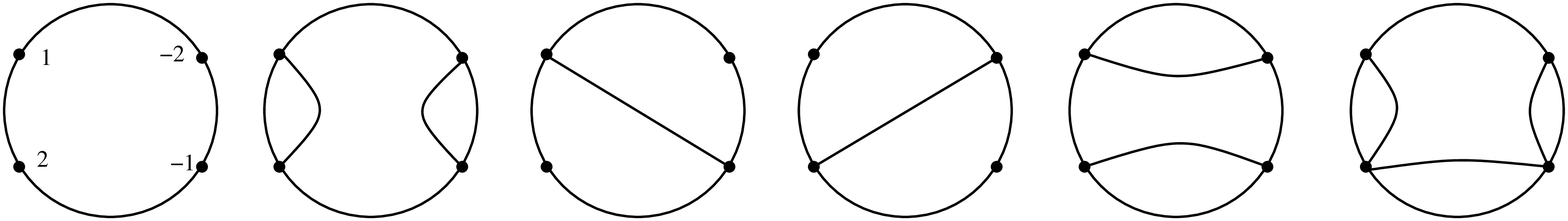}
\caption{The pictorial representation of $\pb2$.}\label{noncrB2}
\end{figure}

\begin{figure}[ht]
\centering
\includegraphics[width=\textwidth]{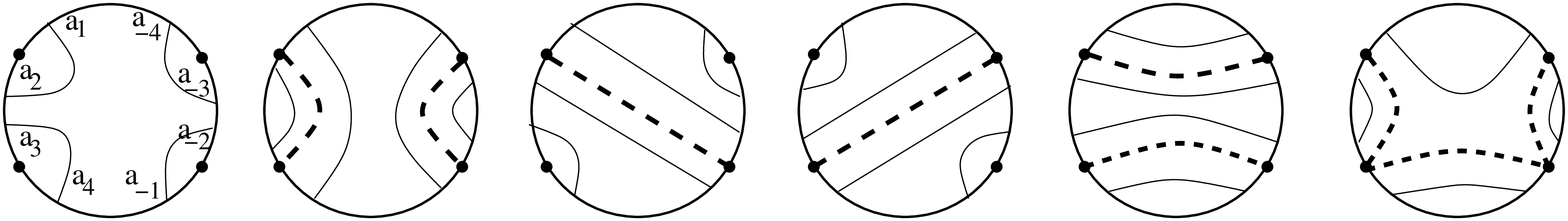}
\caption{The bijection $\bi B:\pb2\to\bb2$.}\label{figBbijection}
\end{figure}

An {\it $n$-partition of type $B$} is a partition $\pi$ of the $2n$ element set 
$\{+1,+2,...,+n,-1,-2,...,-n\}$ into blocks with the property that 
for any block $K$ of $\pi$, its opposite $-K$ is also a block of $\pi$, 
and that there is at most one invariant block (called the 
zero block\footnote{For topologists the term invariant block is 
 more natural than zero block so we use this names interchangeably 
in the paper.}) for 
which $K=-K$. Since all non-zero blocks occur in pairs $\pm K$ one defines 
$\nzbk(\pi)$ as half of the number of all non-zero blocks. To represent $\pi$ 
pictorially, we place the numbers $+1,+2,...,+n,-1,-2,...,-n$ anti-clockwise
around the boundary circle of a disk 
and draw a connection chord in the disk between 
two numbers $i<j$\footnote{Here we use the order 
$+1<+2<\cdots<n<-1<-2<\cdots<-n$.} if they are 
in the same block of $\pi$ and there is no $k$ in 
the same block with $i<k<j$. Then $\pi$ is said to be
{\em non-crossing} if all connection chords can be drawn 
without crossing each other.\footnote{The non-crossing condition 
forces a partition to
have at most one zero block.}
Denote the set of all non-crossing $n$-partitions of type $B$ by $\pb n$.
The set $\pb2$ is illustrated in Fig.~\ref{noncrB2}.
On the other hand if $n=2m$ is even, we have bipartitions of $4m$ points 
of type $B$, those $2m$-partitions of type $B$ with every block 
containing exactly 2 numbers. 
Denote the set of all non-crossing $2m$-bipartitions of type $B$ by $\bb m$.
Similar to the type $A$ case we have a bijection $\bi B:\pb n\to\bb n$ realized by considering the boundary arcs 
of a regular neighborhood of the connection 
chords (see Fig. \ref{figBbijection}).

\subsection{The matrices}
For any $n$-partitions $\pi$ and $\pi'$ of type $A$ (resp. $B$), 
denote by $\pi\vee\pi'$
the finest $n$-partition (not necessarily non-crossing) of 
type $A$ (resp. $B$) that is coarser than both $\pi$ and $\pi'$.
The {\em matrix of chromatic joins} of type $A$ and $B$ are respectively:
$$
(J^A_n(\d))_{\pi,\pi'\in\pa n} = \d^{\bk(\pi\vee\pi')}
\quad \text{and}\quad 
(J^B_n(\d))_{\pi,\pi'\in\pb n} = \d^{\nzbk(\pi\vee\pi')}.
$$ 
For any $2n$-bipartitions $\pi$ and $\pi'$ of type $A$ (resp. $B$), 
one can glue them along the boundary circles respecting the labels. 
The result, denoted $\pi\vee\pi'$, is a collection of disjoint circles 
on a 2-sphere. 
The {\em Gram matrix of Temperley-Lieb algebra} of type $A$ and $B$ 
are respectively\footnote{The matrix $G^A_n(\d)$ was first used by H. Morton 
and P. Traczyk to find a basis of the Kauffman bracket skein module 
of a tangle \cite{MT}, and played an important role in Lickorish's approach 
to Witten-Reshetikhin-Turaev invariants of 3-manifolds \cite{Li}. 
The matrix $G^B_n(1,\d)$ was first considered by Rodica Simion in 1998; 
compare \cite{Sch}.}:
$$
(G^A_n(\d))_{\pi,\pi'\in\ba n} = \d^{c(\pi\vee\pi')} 
\quad \text{and}\quad 
(G^B_n(\a,\d))_{\pi,\pi'\in\bb n} = 
\a^{c_0(\pi\vee\pi')}\d^{c_d(\pi\vee\pi')},
$$
where $c(\pi\vee\pi')$ is the number of circles, $c_0(\pi\vee\pi')$ 
is the number of zero (i.e. invariant) circles $C$ with $C=-C$, and
$c_d(\pi\vee\pi')$ is the number of pairs of non-zero 
circles $C, -C$ with $C\ne -C$ in $\pi\vee\pi'$.

\section{Proof of Theorems \ref{A} and \ref{B}}

\begin{proof}[Proof of Theorem \ref{A}]
For $\pi_i\in\pa n$ let $b_i:=\bi A(\pi_i)\in\ba n$. Since
$c(b_i \vee b_j)$ is also equal to the number of boundary 
components of the regular neighborhood of the pictorial representation of 
$\pi_i \vee \pi_j$ we have
$2\,\bk(\pi_i \vee \pi_j) = c(b_i \vee b_j) + \bk(\pi_i) + \bk(\pi_j) -n =
\bk(\pi_i) - \frac{n}{2} + c(b_i \vee b_j) + \bk(\pi_j) - \frac{n}{2}$. 
The formula can be obtained from  the expression for the Euler 
characteristic of a plane graph: \ 
Let $G_{\pi}$ be a graph 
corresponding to the non-crossing partition $\pi$. 
$G_{\pi}$ is a forest of $n$ vertices 
and $n-\bk(\pi)$ edges. Similarly, let $G_{\pi_i \vee \pi_j}$ be the graph 
corresponding to $\pi_i \vee \pi_j$. We should stress that $\pi_i \vee \pi_j$
does not have to be a noncrossing partition and that the graph 
$G_{\pi_i \vee \pi_j}$ is a plane graph obtained by putting $G_{\pi_i}$ 
inside a disk and $G_{\pi_j}$ outside the disk with $G_{\pi_i} \cap G_{\pi_j}$
composed of the $n$ points on the unit circle\\ (e.g.: 
\parbox{2.5cm}{\includegraphics[height=0.7cm]{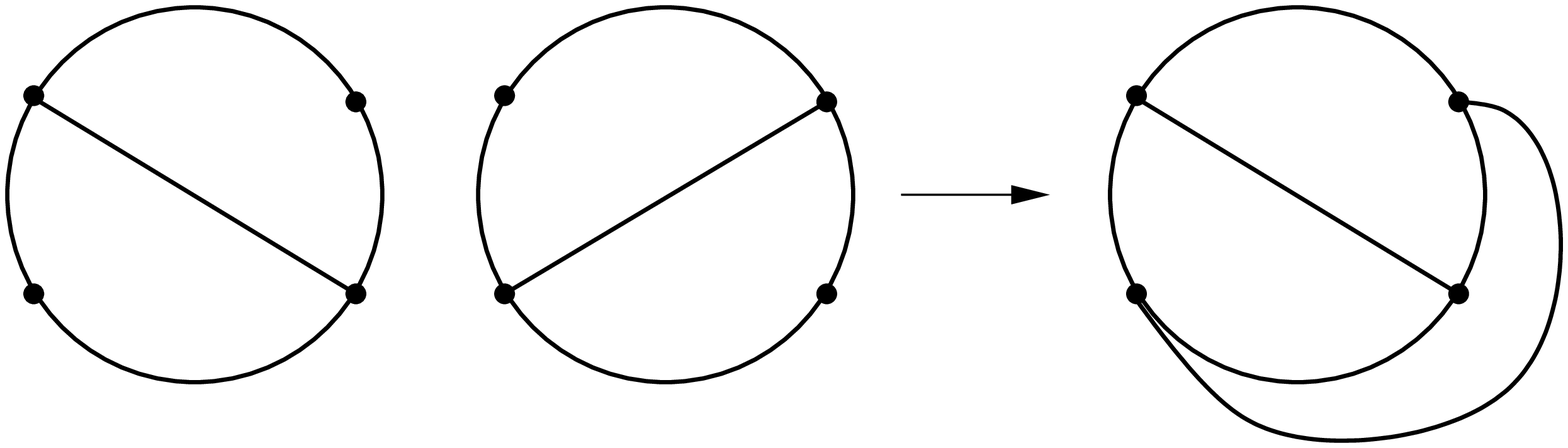}}, or 
\parbox{2.3cm}{\includegraphics[height=0.7cm]{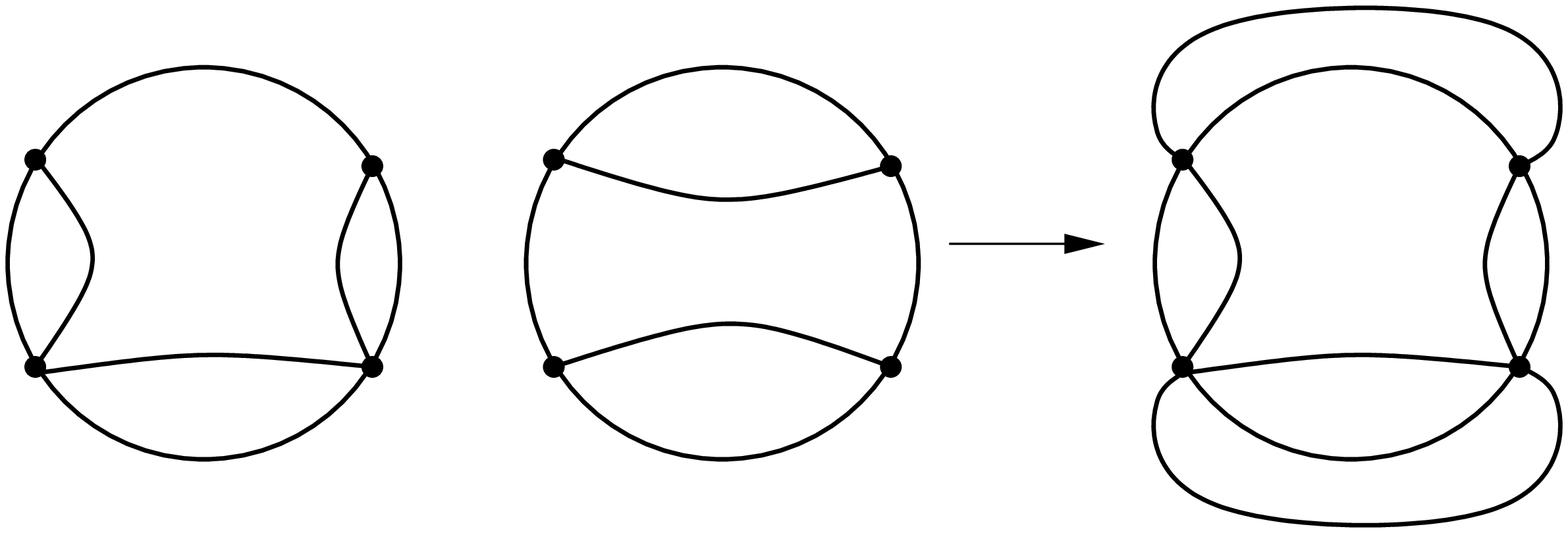}}).\\ 
By construction, $G_{\pi_i \vee \pi_j}$ is a plane graph of $n$ vertices 
and $\bk(\pi_i \vee \pi_j)$ components. It has 
$E(G_{\pi_i \vee \pi_j}) = E(G_{\pi_i})+E(G_{\pi_j})= 
n-\bk(\pi_i) + n-\bk(\pi_j)$ edges. Furthermore, if we embed 
$G_{\pi_i \vee \pi_j}$ in a disjoint union 
of $\bk(\pi_i \vee \pi_j)$ 2-spheres 
(each component of $G_{\pi_i \vee \pi_j}$ in a different sphere) we 
can identify $c(b_i \vee b_j)$ with the number of regions of the 
embedded graph. The Euler characteristic is on the one hand equal to 
$2\bk(\pi_i \vee \pi_j)$ and on the other hand equal to 
$n - E(G_{\pi_i \vee \pi_j}) + 
c(b_i \vee b_j)= c(b_i \vee b_j) + \bk(\pi_i) + \bk(\pi_j)-n$, as 
needed.\\
Theorem \ref{A} follows directly from the formula.
\end{proof}

\begin{proof}[Proof of Theorem \ref{B}]
For $\pi_i\in\pb n$ let $b_i:=\bi B(\pi_i)\in\bb n$ and $\bk_0(\pi_i)$ be 
the number of zero blocks of $\pi_i$.
Recall that 
$c_0(b_i\vee b_j)$ is the 
number of zero (i.e. invariant) components of $b_i \vee b_j$.
As in the type $A$ case we have 
$2\,\bk(\pi_i \vee \pi_j) = c(b_i \vee b_j) 
+ \bk(\pi_i) + \bk(\pi_j) -2n$. 
Furthermore, we have (see Lemma \ref{$bk_0$ Lemma}):
$$
2\,\bk_0(\pi_i \vee \pi_j) = c_0(b_i \vee b_j) + 
\bk_0(\pi_i) + \bk_0(\pi_j).
$$
(Notice that $\bk_0(\pi_i \vee \pi_j)$ 
can be equal to $2$, $1$, or $0$.) From these we conclude that:\\
$2\,\nzbk(\pi_i \vee \pi_j) = c_d(b_i \vee b_j) +
\nzbk(\pi_i) + \nzbk(\pi_j) - n$. Thus Theorem \ref{B} follows.
\end{proof}

\begin{lemma}\label{$bk_0$ Lemma}
The zero blocks and zero components 
satisfy the following identity:
$$
2\,\bk_0(\pi_i \vee \pi_j) = c_0(b_i \vee b_j) +
\bk_0(\pi_i) + \bk_0(\pi_j),
$$
where $b_i = \bi B(\pi_i)$ and $b_j = \bi B(\pi_j)$.
\end{lemma}

\begin{proof}
The lemma reflects the basic properties of a 2-sphere with an 
involution fixing two points and its compact invariant submanifolds.

To demonstrate the formula we consider all cases of 
blocks of $\pi_i$, $\pi_j$ and $\pi_i \vee \pi_j$ divided into four classes:\\
(1) If $K$ is a non-zero (i.e. non-invariant) block of $\pi_i \vee \pi_j$ 
then all its constituent blocks in $\pi_i$ and $\pi_j$ are non-zero blocks 
and the boundary components of a regular neighborhood 
of the geometric realization 
of $K$ (denoted $\partial K$), that is circles in $b_i \vee b_j$ 
corresponding to $K$, are non-invariant (non-zero) curves, i.e. not in 
$c_0(b_i \vee b_j)$. \\
(2) If $K$ is a zero block of $\pi_i \vee \pi_j$ but all 
its constituent blocks in $\pi_i$ and $\pi_j$ are non-zero blocks, then 
exactly two components of $\partial K$ are invariant curves.\\
(3) If $K$ is a zero block of $\pi_i \vee \pi_j$ and exactly one 
constituent block is a zero block then exactly one component of 
$\partial K$ is an invariant curve.\\
(4) If $K$ is a zero block of $\pi_i \vee \pi_j$ and exactly two 
constituent blocks are zero-blocks (necessarily one in $\pi_i$ and one in 
$\pi_j$) then no component of $\partial K$ is an invariant curve.

These conditions taken together prove the formula in Lemma \ref{$bk_0$ Lemma}.
\end{proof}

\section{Corollaries}
Theorem B and the results of \cite{MS,CP} allow us to answer 
Problems 1 and 2 of \cite{Sch} about a formula for the determinant of the 
type-$B$ matrix of chromatic joins:
\begin{corollary} 
$$\det(J^B_n(\d^2))= 
\prod_{i=1}^n \left( T_i(\d)^2 - 1 \right)^{\binom{2n}{n-i}}
$$
where $T_i(\d)$ is the Chebyshev polynomial of the first kind: 
$$
T_0 = 2, \qquad T_1 = \d, \qquad T_i = \d\, T_{i-1} - T_{i-2}.
$$
\end{corollary}
\addtocounter{footnote}{1}

The matrix $J^B_n(\d)$ can be generalized to a matrix of two variables 
as follows:
$$
(J^B_n(\a,\d))_{\pi,\pi'\in\pb n} 
= \a^{\bk_0(\pi \vee \pi')} \d^{\nzbk(\pi \vee \pi')}.$$
It follows from Lemma \ref{$bk_0$ Lemma} that
$$J^B_n(\a^2,\d^2)=P^B_n(\a,\d)G^B_n(\a,\d)P^B_n(\a,\d),$$ 
where $P_n^B(\a,\d)=(p_{ij})$ is a diagonal matrix
with $p_{ii}(\a,\d)= \a^{\bk_0(\pi_i)}\d^{\nzbk(\pi_i)-n/2}$.
Furthermore, $\det P_n^B(\a,\d) = 
\a^{\frac{1}{2}\binom{2n}{n}}$\footnote{This follows from
Proposition 3 of \cite{Rei}, which asserts that there exists 
a fixed-point free involution $\gamma$ on $\pb n$ such that
$\bk_0(\pi)+\bk_0(\gamma(\pi)) = 1$ 
and $\nzbk(\pi)+\nzbk(\gamma(\pi)) = n$. 
} 
and thus we have: 

\begin{corollary}\label{cor}
$$
\det J^B_n(\a^2,\d^2)= \a^{\binom{2n}n} G^B_n(\a,\d) =
\a^{\binom{2n}n}
\prod_{i=1}^n \left( T_i(\d)^2 - \a^2 \right)^{\binom{2n}{n-i}}.
$$
\end{corollary}

\begin{figure}[ht]
\psfrag{a}{$a$}
\psfrag{a'}{$a'$}
\psfrag{S}{\hspace{-.3ex}$S$}
\psfrag{S'}{\hspace{-.5ex}$S'$}
        \centering
                \includegraphics[height=3.9cm]{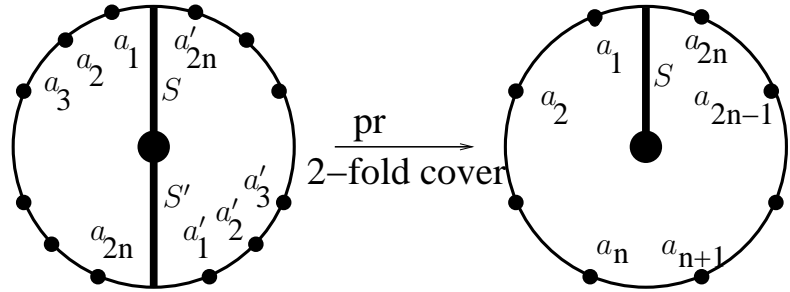}
\caption{The double branch cover $pr: D^2 \to D^2$
with the ``cutting" arc $S$}\label{figBR}
\end{figure}

\begin{rmk} {\em Consider the Gram matrix of type $B$ based
on non-crossing connections in an annulus.\footnote{This interpretation 
of the Gram matrix of type B Temperley-Lieb algebra 
is mentioned in \cite{Sch} as an annular skein matrix 
and utilized in \cite{MS} and \cite{CP}.} This matrix is the same
as the one considered before in Theorem \ref{B} via the branched
cover described in Fig. \ref{figBR}.}
\end{rmk}


\vspace{5ex}
\begin{minipage}[b]{0.45\linewidth}

Department of Mathematics\newline\indent
Winston-Salem State University\newline\indent
Winston Salem, NC 27110, USA\newline\indent
chenqi@wssu.edu

\end{minipage}
\hspace{3ex}
\begin{minipage}[b]{0.45\linewidth}

Department of Mathematics\newline\indent
The George Washington University\newline\indent
Washington, DC 20052, USA\newline\indent
przytyck@gwu.edu

\end{minipage}

\end{document}